\documentclass{ifacconf}

\usepackage{graphicx}      
\usepackage{natbib}        
\usepackage{bm}
\usepackage{amsmath} 
\usepackage{amssymb}  
\usepackage{float}
\usepackage{eqnarray}

\usepackage{algorithmic}
\usepackage{algorithm}
\usepackage{xcolor}

\newtheorem{theorem}{Theorem}

\newtheorem{assumption}{Assumption}

\begin{document}
\begin{frontmatter}

\title{ Preference-based optimization from noisy pairwise comparisons} 

\thanks[footnoteinfo]{*This work was supported in part by Swedish Research Council Distinguished Professor Grant 2017-01078, Knut and Alice Wallenberg Foundation Wallenberg Scholar Grant, and Swedish Strategic Research Foundation SUCCESS Grant FUS21-0026.}

\author[Author_KTH]{Siyi Wang} 
\author[Author_KTH]{Zifan Wang} 
\author[Author_KTH]{Karl Henrik Johanssson}

\address[Author_KTH]{KTH Royal Institute of Technology, Stockholm, Sweden, \\(e-mail: siyiw,  zifanw, kallej@kth.se)}

\begin{abstract}                
In interactive systems, feedback is often provided in the form of preference between queried options rather than precise scores, which motivates optimization methods to learn from such comparisons.
In this work, we propose a preference-based optimization algorithm that relies on noisy two-point comparisons. At each iteration, the algorithm employs a uniform-sphere perturbation to generate a perturbed action and queries the resulting loss comparison to estimate a descent direction. We demonstrate that, under standard smoothness and bounded variance assumptions, the algorithm converges to a stationary point when the smoothing and step size parameters are properly chosen. Numerical experiments on an LQG system demonstrate the effectiveness of the preference-based optimization algorithm with comparison feedback.
\end{abstract}

\end{frontmatter}
\section{Introduction}

Preference-based optimization refers to a learning framework where decisions are improved by comparative judgment or ranking feedback. In many interactive systems, feedback appears in the form of preferences or rankings rather than absolute scores.
For instance, in human-in-the-loop systems, it may be difficult for a human to assign an exact numerical reward to an outcome, but much easier to say “option A is better than option B.” 
In personalized recommender systems, a user’s clicks or choices reveal only relative satisfaction with recommended items, rather than a calibrated rating.  
Such observations motivate the development of optimization algorithms that learn from ordinal comparisons.
Over the past decade, optimization utilizing preference-based and ranking feedback has been investigated across several areas, including dueling bandits \cite{yue2009interactively,kumagai2017regret}, Bayesian optimization \cite{astudillo2020multi}, and deep reinforcement learning \cite{christiano2017deep}.

Optimizing parameters to minimize an objective is a fundamental task across learning and control systems.
When the loss function is available in a closed form, gradient descent provides a principled update rule. However, in many real-world settings, the explicit form of the loss function and, consequently, its gradient are not accessible. This limitation motivates alternative approaches that estimate gradients using indirect feedback, including zeroth-order optimization, which assumes access to function evaluations, and preference-based optimization, which is based on a comparison oracle. Among these methods, the zeroth-order optimization has been extensively explored and supported with a performance guarantee \cite{flaxman2004online,shamir2017optimal}. Existing zeroth-order methods include one-point \cite{flaxman2004online} and multi-point schemes \cite{shamir2017optimal,zhang2022new}; and use perturbations such as Gaussian smoothing \cite{nesterov2017random}, uniform spherical directions \cite{flaxman2004online}, and coordinate-wise perturbations \cite{balasubramanian2018zeroth}.
When only the preference feedback is available, preference-based optimization infers the gradient descent direction from comparative judgments between options. Such optimization method finds application across various fields, including reinforcement learning \cite{ji2024reinforcement}, LLM preference alignment \cite{chen2025compo}, and control systems \cite{wang2025human}. 

This work focuses on developing a gradient descent method with a comparison oracle. 
Within this domain, existing literature includes \cite{kumagai2017regret,ding2018preference,cai2022one,zhang2024comparisons,tang2024zerothorder,lobanov2024acceleration}.
For instance, \cite{cai2022one} considers a pairwise comparison oracle and constructs a robust estimator of the normalized gradient via one-bit compressive sensing. Their approach requires sparsity of the true gradient and convexity of the objective to make the problem tractable. However, confining the objective to convex cases can be unrealistic in human preference applications \cite{tang2024zerothorder}. To address this limitation, \cite{tang2024zerothorder} introduced the ZO-RankSGD algorithm that performs stochastic gradient descent steps using only ranking information and proves convergence to a stationary point under the smooth function assumption.
Additionally, 
\cite{zhang2024comparisons} investigates algorithms via comparison-only access under different loss function settings, including convex and nonconvex smooth functions, achieving convergence rates that match the optimal guarantees for zeroth-order optimization.
Furthermore, 
\cite{lobanov2024acceleration} investigates the accelerated optimization algorithms that rely solely on an order oracle in various settings, which shows that acceleration is possible even with ranking feedback. 
Some work, such as \cite{tang2024zerothorder}, focuses on deterministic objective functions, whereas practical systems inherently involve stochastic cost evaluations. For example, in interactive or human-in-the-loop environments, the evaluation of an action may vary due to human subjectivity or fluctuating environmental conditions. As a result, the pairwise comparisons derived from such evaluations are also stochastic.  
Based on this, this work aims to address optimization under stochastic objectives with noisy pairwise comparisons.

In this paper, we propose an optimization algorithm using feedback through noisy two-point comparisons. At each iteration, the algorithm samples a random direction on the unit sphere and performs a small perturbation in that direction. It then queries the comparison oracle on the current decision and the perturbed point to obtain a binary feedback indicating which one attains a better objective value. Using this feedback, we construct a stochastic gradient estimator to perform gradient descent updates.
We analyze the convergence in terms of the 
gradient norm under standard smoothness and bounded-variance assumptions. The theoretical analysis demonstrates that the algorithm converges to a first-order stationary point under an appropriately chosen step size and smoothing parameter.
Finally, simulations on an LQG system demonstrate the efficacy of gradient descent with the comparison oracle, showing that the control gain converges to the theoretical optimum.

The remainder of this paper is structured as follows: Section \ref{sec:preliminaries} introduces the preliminaries and problem statement. Section \ref{sec:main result} presents the learning algorithm with noisy pairwise comparisons and analyzes its performance. Section \ref{sec:simulation} demonstrates the efficacy of the algorithm by numerical simulations. Section \ref{sec:conclusion} draws conclusions. 
 
\textbf{Notations:}   
Let $\|\cdot\|$ denote the $l_2$ norm.  
Let the notation $\mathcal{O}$ hide the constants and $\tilde{\mathcal{O}}$ hide constants and polylogarithmic factors of the total number of iterations $T$, respectively. 
Let $\mathbb{R}$
 denote the set of real numbers. 
Let $\mathbb{R}^n$ denote the  $n$-dimensional real vector space.
Let $\mathbb{E}$ denote the expectation operator.

\section{Problem statement}\label{sec:preliminaries}

Consider the cost function $f(x,\xi) : \mathbb{R}^d \times \Xi
\rightarrow \mathbb{R} $, where $x\in \mathbb{R}^d$ denotes the decision variable and $\xi \in \Xi$ denotes the random noise with $\Xi$ being its support. Let the random variable $\xi$ follow the distribution $\mathcal{D}$. 
Consider the optimization problem:
\begin{equation}\label{eq:opt}
    \min_{x \in \mathcal{X}}\mathbb{E}_{\xi \sim\mathcal{D}}[  f(x,\xi)].
\end{equation}
For two candidate decisions $x_1,x_2 \in \mathbb{R}^d$, let $f(x_1,\xi)$, $f(x_2,\zeta)$  be the evaluation result, where $\xi, \zeta \sim \mathcal{D}$ represents the randomness in evaluation. 
Assume that the learner cannot directly access function evaluations. Instead, it only observes the relative ordering between two noisy outcomes:
\begin{align}\label{eq:comparison signal}
{\mathrm{sgn}}\big(f(x_1, \xi)-f(x_2,\zeta) \big) \in \{1,-1\} .
\end{align}
Notably, we assume that each comparison is based on two stochastic evaluations. This assumption reflects practical scenarios in interactive or human-in-the-loop systems, where feedback is subject to human variability and fluctuating environmental conditions.

We impose the following assumption, which is standard in the optimization literature \cite{tang2024zerothorder}. 
\begin{assumption}\label{assumption:cost function}
    The function $f$ satisfies: 1) $f$ is twice continuously differentiable; 2) $f$ is $L$-smooth in $x$, i.e., $\|\nabla^2 f(x,\xi)\| \le L$ for all $\xi \in \Xi$, where $\nabla^2 f(x,\xi)$ is the Hessian matrix of $f$ at $x$. 
\end{assumption}
Preference-based optimization algorithms rely on a randomized gradient estimator, which is typically constructed via random finite function queries. 
To ensure convergence and stability toward a stationary point, we provide the following assumption to bound the gradient variance and the function-value variance.  
\begin{assumption}\label{assumption:bounded variance}
Let $F = \mathbb{E}_{\xi}[f(x,\xi)]$. For all $x \in \mathcal{X}$, there exist constants $\sigma,
\sigma_f>0$ such that 
 \begin{align*}
    \mathbb{E}_{\xi}\big[\|\nabla f(x,\xi)-\nabla F(x)\|^2\big] &\le \sigma^2,  \\
    \mathbb{E}_{\xi}\big[   \big(f(x,\xi)\big)^2\big]-  \big(F(x)\big)^2 &\le \sigma_f^2. 
 \end{align*} 
\end{assumption}
The constants 
$\sigma$ and $\sigma_f$ quantify the noise levels in the stochastic gradients and the stochastic function evaluations, respectively.
Similar assumptions can be found in \cite{spall2002multivariate,ghadimi2013stochastic,nesterov2017random}.  

Our goal is to develop a preference-based optimization algorithm that solves \eqref{eq:opt} using the noisy pairwise comparison feedback \eqref{eq:comparison signal}.

\section{Main results}\label{sec:main result}
In this section, we present the optimization algorithm with noisy pairwise comparisons and analyze its convergence properties.
\subsection{Algorithm}
We first present the optimization algorithm with a comparison oracle.   
Given a point $x_t\in \mathcal{X}$, define the perturbed action by 
$\hat{x}_t=x_t+\delta u_t$, where $u_t$ is the direction vector sampled from a unit sphere $\mathbb{S}^d \in \mathbb{R}^d$ at iteration $t$, and $\delta$ is the perturbation radius, also known as the smoothing parameter. Assume that the perturbation vector $u_t$ is independent of the random sample $\xi_t,\zeta_t$, for all $t$.  
Then, we construct the gradient estimate as 
\begin{align}\label{eq:gradient estimate}
g_t = \frac{d}{u_t}{\mathrm{sgn}}\big(f(x_t+\delta u_t,\xi_t)-f(x_t,\zeta_t) \big)u_t,
\end{align}
where $\xi_t$ and $\zeta_t$ follow the distribution of $\mathcal{D}$.
The gradient estimator \eqref{eq:gradient estimate} is inspired by the standard two-point zeroth-order optimization method \cite{shamir2017optimal} as well as recent developments in optimization with ranking information \cite{tang2024zerothorder}. 

Furthermore, the gradient descent proceeds as 
\begin{equation}\label{eq:gradient descent}
x_{t+1} = x_t - \eta g_t
\end{equation}
with the initial value $x_0$ and the step size $\eta$. The learning algorithm with the noisy pairwise comparison feedback is summarized in Algorithm~\ref{alg:algorithm}. 
\begin{algorithm} 
\caption{Noisy pairwise comparison SGD} \label{alg:algorithm}
\begin{algorithmic}
    \REQUIRE Initial value $x_0$, iteration horizon $T$, smoothing parameter $\delta$, step size $\eta$.
    \FOR{$ {\rm{iteration}} \;t = 0,\dots, T$} 
    \STATE  Sample $u_{t} \in \mathbb{S}^{d}$ 
    \STATE  Play $x_t$, $ x_{t}+\delta u_{t} $  and obtain comparison feedback ${\rm sgn}\big(f_t(x_{t}+\delta u_{t} ,\xi_t)-f_t(x_t,\zeta_t)\big)$
    \STATE Construct gradient estimate $g_{t}$, as in \eqref{eq:gradient estimate}
    \STATE Update decision: $x_{t+1}=  x_{t} - \eta  g_{t}$
    \ENDFOR
\end{algorithmic}
\end{algorithm}

\subsection{Convergence analysis} 
This section analyzes the convergence performance of Algorithm~\ref{alg:algorithm}.  
In the context of smooth, nonconvex optimization, the primary objective is to prove the algorithm's convergence to an approximate stationary point. Since the iterative update is based on a noisy, randomized gradient estimator, the state $x_t$ is a random variable dependent on the history of all prior random samples. 
Define by $\mathcal{Z}_t=\{u_0, \xi_0,\zeta_0,\dots,u_t, \xi_t,\zeta_t\}$ a random vector composed by variables $z_t=\{u_t,\xi_t,\zeta_t\}$  attached to each iteration of the Algorithm. Moreover, denote $G_t:= \mathbb{E}_{\mathcal{Z}_{t-1}}\big[\|\nabla F(x_t)\|\big]$, for all $t$. 
Therefore, we analyze the property of the average conditional expected gradient norm, i.e., $\frac{1}{T}\sum_{t=1}^T G_t$, in the following theorem. 
\begin{theorem}
Let Assumptions~\ref{assumption:cost function} and \ref{assumption:bounded variance} hold. Select the step size $\eta = T^{-\frac{1}{2}}$ and the smooth parameter as a positive constant. The convergence rate for the average expected gradient norm under Algorithm~\ref{alg:algorithm} satisfies 
\begin{equation*}
\frac{1}{T}\sum_{t=1}^T G_t  = \tilde{\mathcal{O}}(T^{-\frac{1}{2}} + \sigma +\sigma_f). 
\end{equation*}
\end{theorem}
\textit{Proof.}
Since the function $f$ is $L$-smooth in $x$, we have that 
\begin{align}\label{eq:smooth inequality}
F(x_{t+1}) &= F(x_t-\eta g_t) \nonumber \\ 
&\le F(x_t) - \eta \langle \nabla F(x_t),g_t\rangle + \frac{L\eta^2}{2}\|g_t\|^2,
\end{align}
where $F(x)=\mathbb{E}_{\xi}[f(x,\xi)]$. 
Applying the expectation operator to \eqref{eq:smooth inequality} yields  
\begin{align}\label{eq:exp f}
&\mathbb{E}_{z_t}[F(x_{t+1})] \nonumber\\
\le &  \mathop F(x_{t})-\eta  \mathbb{E}_{z_t}[\langle \nabla F(x_{t}), g_{t}\rangle]+\frac{L \eta^{2}}{2}  \mathbb{E}_{z_t}[\|g_{t}\|^{2}],   
\end{align}
where $\nabla F(x)$ denote the derivative of $F$ with respect to $x$. 
We next bound the second term of the right-hand side of \eqref{eq:exp f}. 
It follows that 
\begin{align}\label{eq:int_df_g_1}
&\mathbb{E}_{z_t}\big[\langle \nabla F(x_t),g_t\rangle \big] \nonumber \\
=& \mathbb{E}_{z_t}\big[\langle \nabla F(x_t),u_t \rangle \frac{d}{\delta}{\rm sgn}\big(f(x_t+\delta u_t,\xi_t)-f(x_t,\zeta_t) \big)\big] \nonumber \\
=&  \frac{d}{\delta} \mathbb{E}_{z_t}\Big[\langle \nabla F(x_{t}), u_{t}\rangle \operatorname{sgn}\Big(\langle\nabla f(x_{t},\xi_t), \delta u_{t}\rangle \nonumber \\ 
&+\frac{\delta^{2}}{2} u_{t}^{\top} \nabla^{2} f(x_{t},\xi_t)u_{t}+f(x_t,\xi_t) -f(x_t,\zeta_t) \Big)\Big] ,
\end{align}
where the first equality is from the definition of $g_t$, as in \eqref{eq:gradient estimate}, and the second equality is from the second-order Taylor expansion with Cauchy remainders: 
\begin{align*} 
&f(x_t + \delta u_t,\xi) \nonumber \\
=& f(x_t,\xi) + \delta \langle \nabla f(x_t,\xi),u_t\rangle + \frac{\delta^2}{2}u_t^\top\nabla^2f(x_t,\xi)u_t 
\end{align*}
for all $\xi \in \Xi$. 
Furthermore, denote  $A=\langle\nabla f(x_{t},\xi_t), \delta u_{t}\rangle$ and $B = \frac{\delta^2}{2} u_{t}^{\top} \nabla^{2} f_{t}(x_{t},\xi_t) u_{t} + f_{t}(x_{t},\xi_t)-f_{t}(x_{t},\zeta_t)$. Then, \eqref{eq:int_df_g_1} is further written as 
\begin{align}\label{eq:int_df_g_2} 
    & \quad \mathbb{E}_{z_t}\big[\langle \nabla F(x_t),g_t\rangle \big] \nonumber \\
    &= \frac{d}{\delta^2} \mathbb{E}_{z_t}\Big[ \big(\langle\nabla F(x_t),\delta u_t \rangle - A \big)\operatorname{sgn}(A+B) \nonumber \\
    &\quad + \big|A+B\big| - B\operatorname{sgn}(A+B)\Big] \nonumber \\
    &\ge \frac{d}{\delta^2} \mathbb{E}_{z_t}\big[(\langle\nabla F(x_t),\delta u_t \rangle - A )\operatorname{sgn}(A+B) + |A|-2|B| \big],
\end{align}
where the equality is from $ |A+B|= (A+B){\rm sgn}(A+B) $,  and the inequality is from  $|A+B| \ge  |A|-|B|$ and $B{\rm sgn}(A+B)\ge -|B|$. 

We next bound the terms in the right-hand side of \eqref{eq:int_df_g_2}, i.e., $ \mathbb{E}_{z_t}\Big[\big(\langle\nabla F(x_t),\delta u_t \rangle -A \big){\rm sgn}(A+B)  \Big]$, $\mathbb{E}_{z_t}[|A|]$ and $-2\mathbb{E}_{z_t}[|B|]$, respectively. 
For the first term, we have that 
\begin{align}\label{eq:func derv cov bound}
&\mathbb{E}_{z_t}\Big[\big(\langle\nabla F(x_t)-\nabla f(x_t,\xi_t),\delta u_t \rangle  \big){\rm sgn}(A+B)  \Big]   \nonumber\\
\ge & -\mathbb{E}_{u_t,\xi_t}\Big[ \big|\langle\nabla F(x_t)-\nabla f(x_t,\xi_t),\delta u_t \rangle \big|  \Big]  \nonumber \\ 
\ge& - \delta \mathbb{E}_{u_t,\xi_t}\big[\| \nabla F(x_t)-\nabla f(x_t,\xi_t) \|  \|u_t\|\big] \nonumber \\ 
\ge&  - \delta \sqrt{\mathbb{E}_{\xi_t}\big[\| \nabla F(x_t)-\nabla f(x_t,\xi_t) \|^2 \big]  \mathbb{E}\big[  \|u_t\|^2\big]} \nonumber \\
\ge& -\delta \sigma ,
\end{align}
where the last inequality is from applying Cauchy–Schwarz inequality and the fact that the random vector $u_t$ is independent of $\xi_t$.
For  $\mathbb{E}_{z_t}[|A|]$, we have
\begin{align}
\mathbb{E}_{z_t}[|A|]=& \mathbb{E}_{u_t,\xi_t}\big[|\langle \nabla f(x_{t},\xi_t), \delta u_{t}  \rangle| \big] \nonumber \\ 
=& \mathbb{E}_{\xi_t}\Big[\mathbb{E}_{u_t}\big[ |\langle \nabla f(x_{t},\xi_t), \delta u_{t}  \rangle| \big] \Big| \xi_t  \Big]  \nonumber \\ 
\ge &c_d  \delta \mathbb{E}_{\xi_t}\big[\|\nabla f(x_{t},\xi_t)\|  \big] \nonumber\\
\ge& c_d  \delta \|\mathbb{E}_{\xi_t}[\nabla f(x_t,\xi_t)] \| = c_d\delta \| \nabla F(x_{t})\| ,
\end{align}
where $c_d$ is a positive constant depending on the dimension $d$. The first inequality is from $\|u_t\|=1$ 
and the second inequality is from applying Jensen's inequality. 
Furthermore, $\mathbb{E}_{z_t}[|B|]$ is written as 
\begin{align}\label{eq:E[B]}
&\mathbb{E}_{z_t}[|B|]= \mathbb{E}_{z_t}\Big[\Big|\frac{\delta^2}{2} u_{t}^{\top} \nabla^{2} f_{t}(x_{t},\xi_t) u_{t} + f_{t}(x_{t},\xi_t)-f_{t}(x_{t},\zeta_t)\Big|\Big] \nonumber \\ 
&\le  \mathbb{E}_{z_t}\Big[\Big|\frac{\delta^2}{2} u_{t}^{\top} \nabla^{2} f_{t}(x_{t},\xi_t) u_{t}\Big| + \Big|f_{t}(x_{t},\xi_t)-f_{t}(x_{t},\zeta_t)\Big|\Big] \nonumber \\ 
&\le \frac{\delta^2L}{2} + \mathbb{E}_{\xi_t,\zeta_t}\Big[  \Big|f_{t}(x_{t},\xi_t)-f_{t}(x_{t},\zeta_t)\Big|\Big],
\end{align}
where the second inequality is from $\mathbb{E}_{u_t}[u_t^\top \nabla^2 f_t(x_t)u_t]\le \| \nabla^2 f_t(x_t)\| \mathbb{E}[\|u_t\|^2] \le L$, which holds by Assumption~\ref{assumption:cost function}.
Denote $\mu = \mathbb{E}_{\xi}[f(x_t,\xi)]$. Then, the second term of  \eqref{eq:E[B]} is written as 
\begin{align}
&\mathbb{E}_{\xi_t,\zeta_t}\Big[  \Big|f_{t}(x_{t},\xi_t)-f_{t}(x_{t},\zeta_t)\Big|\Big] \nonumber \\
\le& \mathbb{E}_{\xi_t,\zeta_t}\Big[  \big|f_{t}(x_{t},\xi_t)-\mu \big| + \big|f_{t}(x_{t},\zeta_t)-\mu \big|\Big]\nonumber \\ 
=& 2\mathbb{E}_{\xi_t}\big[  \big|f_{t}(x_{t},\xi_t)-\mu \big|\big] \nonumber\\
\le& 2\sqrt{\mathbb{E}_{\xi_t}[(f_{t}(x_{t},\xi_t)-\mu )^2]} \le  2\sigma_f
\end{align}
where the second inequality is from applying Cauchy–Schwarz inequality and the last inequality is from the bounded function variance, as shown in Assumption~\ref{assumption:bounded variance}. 

Additionally, from the definition of gradient estimate \eqref{eq:gradient estimate}, we have that 
\begin{equation}\label{eq:g^2}
\mathbb{E}_{z_t}[\|g_{t}\|^{2}]=\frac{d^{2}}{\delta^{2}} \mathbb{E}[\|u_{t}\|^{2}]=\frac{d^{2}}{\delta^{2}}   .  
\end{equation}
Substituting \eqref{eq:func derv cov bound}-\eqref{eq:g^2} into \eqref{eq:exp f}, we obtain 
\begin{align}\label{eq:Exp f 2}
\mathbb{E}_{z_t}[F(x_{t+1})] & \le F(x_{t})-\frac{c_d\eta d}{\delta} \|\nabla F(x_{t})\| \nonumber \\
 &+ \eta d L +\frac{\eta d\sigma}{\delta} + \frac{4\eta d \sigma_f}{\delta^2}+\frac{L \eta^{2}d^2}{2 \delta^{2}}.
\end{align}
Taking the expectation of \eqref{eq:Exp f 2} in $\mathcal{Z}_{t-1}$, we obtain 
\begin{align}\label{eq:Exp f 3}
\mathbb{E}_{\mathcal{Z}_t}[F(x_{t+1})]  \le & \mathbb{E}_{\mathcal{Z}_{t-1}}[F(x_{t})]-\frac{c_d\eta d}{\delta} G_t\nonumber \\
& + \eta d L +\frac{\eta d\sigma}{\delta} + \frac{4  \eta d\sigma_f}{\delta^2}+\frac{L \eta^{2}d^2}{2 \delta^{2}}.
\end{align}
Rearranging \eqref{eq:Exp f 3} and summing it over $t=1,\dots,T$ yields  
\begin{align*}
\frac{c_d\eta d }{\delta} \sum_{t=1}^{T}G_t    \le & \sum_{t=1}^{T}\big(\mathbb{E}_{\mathcal{Z}_{t-1}}[F(x_{t})]-\mathbb{E}_{\mathcal{Z}_{t}} [F(x_{t+1})\big] \nonumber 
\\
& + \eta  d L T + \frac{\eta d\sigma T}{\delta} + \frac{4\eta d \sigma_f T}{\delta^2} +\frac{L \eta^{2}d^2T}{2 \delta^{2}} \nonumber \\
\le & \mathbb{E}_{\mathcal{Z}_{0}}[F(x_{1})]- \mathbb{E}_{\mathcal{Z}_{T-1}}[F(x_{T})]+  \eta  d L T \nonumber \\ 
&+ \frac{\eta d\sigma T}{\delta} + \frac{4\eta d \sigma_f T}{\delta^2} +\frac{L \eta^{2}d^2T}{2 \delta^{2}}.
\end{align*}
It follows that
\begin{align*}
    \frac{1}{T} \sum_{t=1}^{T}G_t 
   \le&  \frac{\delta\big(\mathbb{E}_{\mathcal{Z}_0}[F(x_{1})]- \mathbb{E}_{\mathcal{Z}_{T-1}}[F(x_{T})]\big)}{c_d\eta d T} \nonumber \\ 
   &+ \frac{1}{c_d} \Big(L\delta +\sigma  + \frac{4\sigma_f+L \eta d}{ \delta} \Big). 
\end{align*}
Then, by selecting $\delta=T^{0}$, $\eta = T^{-\frac{1}{2}}$, Algorithm~~\ref{alg:algorithm} achieves  $\frac{1}{T} \sum_{t=1}^{T}G_t =\tilde{\mathcal{O}}(T^{-\frac{1}{2}}+\sigma+\sigma_f)$, where $\sigma$ and $\sigma_f$ capture the inherent stochasticity in the gradient and function evaluations.  
This bound indicates that the algorithm attains convergence rate $\tilde{\mathcal{O}}(T^{-\frac{1}{2}})$ for stochastic nonconvex optimization. The average stationarity measure decreases as the number of iterations $T$ increases, up to additional error terms stemming from the noise in gradient estimation.

Additionally, in the noiseless setting, i.e., when considering a deterministic function $f(x)$, we have $\sigma=\sigma_f=0$. Substituting these conditions into the result, we recover the noiseless rate of $\tilde{\mathcal{O}}(T^{-1/2})$ \cite{tang2024zerothorder}.

\section{Simulation Result}\label{sec:simulation}
\begin{figure}
    \centering
    \includegraphics[width=0.96\linewidth]{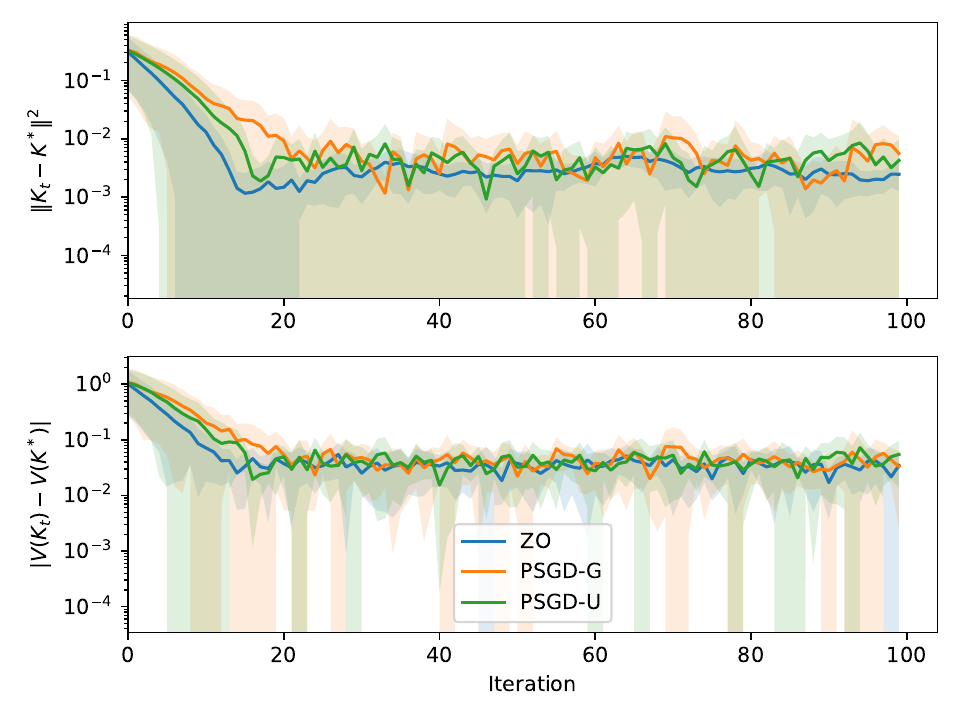}
    \caption{Performance of the algorithms on the system \eqref{eq:plant} with noise variance $\sigma=0.01$.
}
    \label{fig:1DV1}
\end{figure}

\begin{figure}
    \centering    \includegraphics[width=0.96\linewidth]{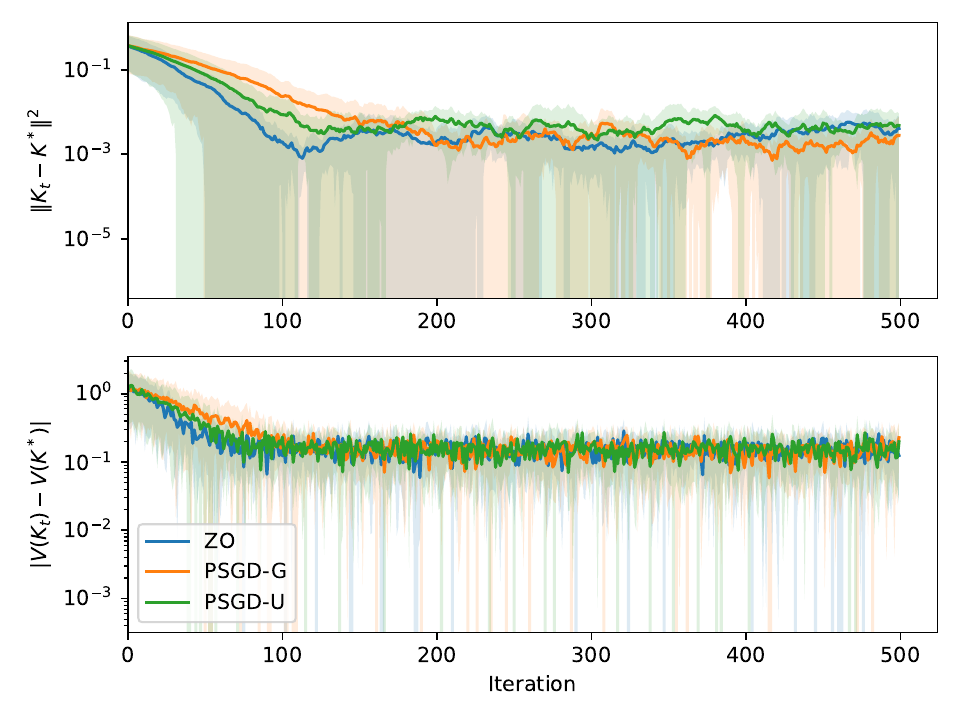}
    \caption{Performance of the algorithms on the system \eqref{eq:plant} with noise variance $\sigma=0.05$.
}
    \label{fig:1DV5}
\end{figure}
The numerical experiments are conducted on a discrete-time LQR system:
\begin{equation}\label{eq:plant}
x_{k+1} = A x_k + B u_k + w_k, \quad u_k = K x_k,    
\end{equation}
where the state $x_k \in \mathbb{R} $ and control $u_k \in \mathbb{R} $. The process noise is set to $w_k \sim \mathcal{N}(0, \sigma^2  )$. The system matrices are chosen as $A=1.1$ and $B=0.1$, respectively.  
The performance is measured by the discounted cumulative cost:
$$V(K) =  \sum_{t=0}^{T} \gamma^t (x_t^\top Q x_t + u_t^\top R u_t) ,$$
where the discount factor is set as $\gamma = 0.7$ and the weighting matrices are $Q= R=1$. The iteration horizon is chosen as $T=50$. 
We first compute the optimal LQR gain $K^*$ analytically via the discrete algebraic Riccati equation. Let $K_t$ denote the policy parameter obtained at iteration $t$. The initial policy for each trial, i.e., $K_0$, is a perturbed version of the optimal policy to ensure a fair starting point:
$K_0 = K^* + \text{Uniform}[0, 1].$  We evaluate the performance of the gradient estimators using: 1) the policy parameter error $\|K_t-K^\ast\|$ and 2) the cost value error  $|V(K_t) - V(K^*)|$.  A smaller policy parameter error indicates that the learning algorithm has driven the policy parameters closer to the known optimal policy matrix $K^\ast$. Moreover, the cost value error measures the quality of the learned policy in terms of the actual objective function, where a smaller value indicates a better-performing control law that minimizes the discounted expected cost. 

Algorithm~\ref{alg:algorithm} uses the preference-based stochastic gradient descent under uniform perturbation (PSGD-U). We compare it with two alternative learning algorithms:
1) ZO: a two-point zeroth-order optimization utilizing bandit feedback \cite{shamir2017optimal}, which uses gradient estimate $\hat{g}_t= \frac{d}{u_t} \big(f(x_t+\delta u_t,\xi_t)-f(x_t,\zeta_t) \big)u_t$ and serves as the benchmark in simulation results; 2) PSGD-G: a preference-based stochastic gradient descent with Gaussian perturbations \cite{tang2024zerothorder}, which uses gradient estimate $\tilde{g}_t=\frac{d}{\tilde{u}_t}{\mathrm{sgn}}\big(f(x_t+\delta \tilde{u}_t,\xi_t)-f(x_t,\zeta_t) \big)\tilde{u}_t,$ where $\tilde{u}_t$ follows a normal distribution.  
The learning rate $\eta$ and the smoothing parameter $\delta$ for each method are carefully tuned to achieve optimal performance. We perform $10$ independent trials for each algorithm and compute their mean and $\pm$ standard deviation. 
The performances of the three algorithms are compared under noise standard deviations $\sigma = \{0.01, 0.05 \}$, which are summarized in Figs.~\ref {fig:1DV1}-\ref{fig:1DV5}.  
Overall, it can be observed that both the policy parameter estimate error and the cost value error converge under three algorithms. The convergence rate slows down as the number of iterations increases. In the preference-based algorithms, though the magnitude of the gradient estimator remains nearly constant, i.e., $\|g_t \| = \frac{d}{\delta}$, its expectation $\mathbb{E}[g_t]$ gradually approaches zero, resulting in slower progress toward the optimum. Moreover, it can be observed that the policy parameter estimation error continues to oscillate even when the cost-value error stabilizes. This occurs because $\mathbb{E}[V(K)]$ is non-convex in $K$, and its gradient remains close to zero in the neighborhood of the optimum. Additionally, the uniform-perturbation algorithm attains a convergence rate comparable to that of the benchmark method, and both converge faster than the algorithm based on Gaussian perturbations. 


\section{CONCLUSIONS}\label{sec:conclusion}
This paper developed a preference-based optimization method that operates with noisy pairwise comparison feedback. We employed a two-point comparison oracle that perturbs the current iterate along a uniformly sampled spherical direction to generate the comparison feedback. Using this feedback, we construct a stochastic gradient estimator for gradient descent updates.
Under smoothness and bounded variance assumptions for stochastic objectives, we establish that the algorithm converges to a first-order stationary point when the smoothing parameter and step size are appropriately chosen. 
Finally, our numerical results demonstrate that the algorithm can effectively optimize stochastic functions using only noisy comparison feedback. These findings underscore the practicality of preference-based optimization with noisy pairwise comparisons, opening up opportunities for broader applications in human-in-the-loop optimization and preference-based learning systems.


\bibliography{ifacconf}          

\appendix
\end{document}